\documentclass[12pt]{amsart}
\usepackage{amssymb,eucal,graphicx}
\usepackage{ifthen}

\usepackage[all,poly,import,color,ps,dvips]{xy}
\def\colorxy(#1){/xycolor{#1 setrgbcolor}def}

\SelectTips{cm}{}

\numberwithin{equation}{section}
1

\newtheorem*{namedtheorem}{\theoremname}
\newcommand{\theoremname}{testing}

\newtheorem{theorem}{Theorem}[section]
\newtheorem{proposition}[theorem]{Proposition}
\newtheorem{proposition-definition}[theorem]
{Proposition-Definition}
\newtheorem{corollary}[theorem]{Corollary}
\newtheorem{lemma}[theorem]{Lemma}
\theoremstyle{definition}
\newtheorem{definition}[theorem]{Definition}

\newtheorem{example}[theorem]{Example}

\newtheorem{remark}[theorem]{Remark}
\theoremstyle{remark}

\def\displaytimes_#1{\mathrel{\mathop{\times}\limits_{#1}}}
\def\displayotimes_#1{\mathrel{\mathop{\bigotimes}\limits_{#1}}}

\newcommand\ext{\operatorname{Ext}}

\newdir{ >}{{}*!/-5pt/@{>}}

\def\p{\mathbb P^2}
\def\n{\mathbb P^{N}}

\def\v{\mathcal{V}_{\d, n}}
\def\d{\delta}
\def\g{\frac{(n-1)(n-2)}{2}}
\def\kn{H^0(\p,\mathcal{O}_{\p}(kH))} 

\def\cp{C^{\prime}}
\def\bc{\begin{center}}
\def\ec{\end{center}}

\def\vuoto{\ \hfill\hbox{\vbox{\hrule\hbox{\vrule
height5pt\kern5pt\vrule height5pt}\hrule}}\par\medskip\rm}

\newcommand{\os}{\mathcal{O}_S}
\newcommand{\ost}{\mathcal{O}_{\tilde{S}}}
\newcommand{\ot}{\mathcal{O}_{\tilde{C}}}

\begin{document}

\title{Geometric $k$-normality of curves and applications}

\author{Alessandro Arsie}
\address{Dipartimento di Matematica\\
Universit\`a di Bologna\\
40126 Bologna\\ Italy}
\email{arsie@dm.unibo.it}
\author{Concettina Galati}
\address{Dipartimento di Matematica\\
Universit\`a degli Studi di Roma Tor Vergata\\
00133 Roma\\ Italy}
\email{galati@axp.mat.uniroma2.it}

\begin{abstract}
The notion of geometric $k$-normality for curves is introduced in complete 
generality and is investigated in the case of nodal and cuspidal curves living 
on several types of surfaces. We discuss and suggest some applications of
this notion to the study of Severi varieties of nodal curves on surfaces of 
general type and on $\mathbb{P}^2$.  
\end{abstract}

\date{November 5, 2003}

%\subjclass[2000]{19E08; 14L30}
\thanks{Authors are partially supported by the University of
Bologna and by the University of Rome Tor Vergata, 
funds for selected research topics.}

%{\small Keywords:geometric $k$-normality, Severi verieties}

\maketitle
%\tableofconten
%%%%%%%%%%%%%%%%%%%%%%%%%%%%%%%%%%
\section{Introduction and preliminaries}
A very classical problem in algebraic geometry is the investigation on  
properties of projectively normal varieties. A closed subscheme 
$X\subset \mathbb{P}^r$ is {\em projectively normal} iff it is normal and 
$H^1(\mathcal{I}_{X/\mathbb{P}^r}(kH))=0$, for all $k\in\mathbb{Z}$ and where 
$H$ denotes the hyperplane section. 
Observe that, when $X$ is smooth, the geometric meaning of the above vanishing 
is that the linear system cut out on $X$ by the hypersurfaces in $\mathbb{P}^r$ 
of degree $k$ is complete. Particular examples of projectively normal subvarieties are 
smooth complete intersections.  For smooth subvarieties, 1-normality (equiv. {\em 
linear normality}), 
means that $X$ cannot be the isomorphic projection of some smooth, non-degenerate variety 
sitting  in a higher dimensional projective space.  

In \cite{cs} it has been introduced a generalization of the concept of linear 
normality, namely the {\em geometric linear normality}. Precisely, in the case of curves, 
we recall the following:
\begin{definition}\label{gln} 
(See \cite{cs}, Def. 3.1). Let $C \subset \mathbb{P}^r$ be any non-degenerate, 
reduced, irreducible curve (possibly singular). 
Then $C$ is said to be geometrically linearly normal if the normalization map $\phi: \tilde{C} \rightarrow C$ 
cannot be factored into a non-degenerate birational morphism  
$\psi: \tilde{C}\rightarrow \Gamma \subset \mathbb{P}^N$, 
followed by a projection $\pi: \mathbb{P}^N \rightarrow \mathbb{P}^r$, $N>r$. 
Equivalently, the linear system $\mathcal{G}(1)$ on $\tilde{C}$, which is the 
pull-back of the linear system cut out on $C$ by the hyperplane sections, is 
complete; i.e. $h^0(\tilde{C},\mathcal{O}_{\tilde{C}}(\phi^{*}(H)))=dim(\mathcal{G}(1))+1=r+1$.
\end{definition}
Let us clarify the sense of the previous definition with the following 
remark. Let $p_a$ be the arithmetic genus of $C$ and $g$ its geometric 
genus. When $p_a$ and $g$ are close enough, then the factorization $\phi=
\pi\circ \psi$ cannot exist, and so $C$ is geometrically linearly normal. 
For example consider a smooth curve $\tilde{C}\subset \mathbb{P}^3$, which lies
on a smooth quadric $Q_2$. Assume that $\tilde{C}$ is of bi-degree $(k,k)$, so that 
its degree is $d=2k$ and $g_{\tilde{C}}=(k-1)^2$. Let us take $d>6$, i.e. $k>3$, and consider the general 
projection $C$ of $\tilde{C}$ to $\mathbb{P}^2$. 
Obviously $C$ is birational to $\tilde{C}$ and we can assume that it has only nodes as singularities, 
since we are using a general projection. 
Thus $p_a(C)=(2k-1)(k-1)$, so that the number $\delta$ of nodes of $C$ is
$\delta=p_a(C)-g_{\tilde{C}}=k(k-1)$. In particular, by using Castelnuovo's 
bound we see that an irreducible nodal curve of degree 8 in $\mathbb{P}^2$ can 
not be the projection of a {\em smooth} curve in $\mathbb{P}^3$, if it has less than 12 nodes. 
Observe that one cannot hope to get this kind of nodal curve in $\mathbb{P}^2$ even using a projection from a higher dimensional 
projective space: indeed by Castelnuovo's theory if we start from a smooth curve $\tilde{C}$ in $\mathbb{P}^r$, instead of $\mathbb{P}^3$, then $g_{\tilde{C}}$
becomes smaller with respect to the degree, so that the number of nodes has to 
increase. Let us remark however that, from Definition \ref{gln}, the above reasoning does not imply that such a curve $C$ in $\mathbb{P}^2$ is geometrically linearly normal; 
instead we have proved that $C$ cannot be the 
projection of a smooth curve, but in principle it could be the projection of 
an already singular curve $\tilde{C}$. Anyhow, this cannot happen by Theorem 2
of \cite{fm}.

On the other hand, if the numer of nodes $\delta$ is big enough, there are 
several examples of irreducible nodal curves on surfaces in $\mathbb{P}^r$ such that $\mathcal{G}(1)$ is not complete (see \cite{cs} for nodal 
curves on surfaces in $\mathbb{P}^3$, \cite{fm} for nodal curves on surfaces in  $\mathbb{P}^r$).
In this framework of ideas, it is natural to try to extend the notion of 
geometric linear normality to that of {\em geometric k-normality}. This is the content of the following:

\begin{definition}\label{gkn}
Let $C\subset \mathbb{P}^r$ be a non-degenerate, reduced, irreducible curve 
(possibly singular) and let $\phi: \tilde{C} \rightarrow C$ be the normalization map. 
Then $C$ is said to be {\em geometrically $k$-normal} (g.k.n., for short) if the pull-back to the 
normalization $\tilde{C}$ of 
the linear system cut out on $C$ by the hypersurfaces of degree $k$ is complete, i.e. of 
dimension $h^0(\mathcal{O}_{\mathbb{P}^r}(k))-h^0(\mathcal{I}_{C/\mathbb{P}^r}(k))-1$. 
As in Definition \ref{gln} the pull-back to $\tilde{C}$ of this linear system 
will be denoted by $\mathcal{G}(k)$.
\end{definition}

\begin{remark}\label{veronese} 
It is immediate to see that $C\subset \mathbb{P}^r$ is 
geometrically k-normal iff  its image $v_k(C)
\subset \mathbb{P}^N$, via the k-fold Veronese (re)embedding $v_k:
C\subset \mathbb{P}^r \hookrightarrow \mathbb{P}^N$, is 
geometrically linearly normal in the sense of Definition \ref{gln}.
\end{remark}

Let us observe that, while the concept of geometrical linear normality does not depend on the postulation of the curve, this is instead the case for geometric 
$k$-normality.

From now on let $S\subset \mathbb{P}^r$ be a smooth, irreducible, non-degenerate k-normal surface. Let $|D|$ be a complete linear system on $S$ whose general element is smooth and irreducible (so called {\em Bertini linear system}) and assume moreover that this general element is $k$-normal. Let $H$ be the hyperplane 
divisor class on $S$ and $K_S$ the canonical divisor class of $S$.
Let $C\in |D|$ on $S$ such that $C$ is a reduced irreducible curve having only 
$\delta$ nodes and cusps as singularities. Let $N$ be the zero dimensional 
subscheme of singularities of $C$. 

The paper is organized as follows. In Section 2 we give necessary and sufficient conditions for such a curve $C$ to be 
geometrically $k$-normal: $C$ is g.$k$.n. iff $N$ imposes independent 
conditions on the linear system $|D+K_S-kH|$ (cfr. Theorem \ref{thm1}). 
Next we find sufficient numerical conditions on $D$, $S$ and $H$ 
which determine an upper-bound $f(H,D,S)$ on 
$\delta$, such that if $\delta<f(H,D,S)$, then $N$ imposes independent 
conditions on $|D+K_S-kH|$. Then, a reduced irreducible curve $C \in |D|$ with $\delta$ nodes and cusps as the only singularities is 
is geometrically k-normal. Precisely, we prove:

\vskip 5pt

\noindent
{\bf Theorem} (cfr. Theorem \ref{thm2} and Corollary \ref{aiuto}) 
{\em Let $k$ be a positive integer and let 
$S\subset \mathbb{P}^{r}$ be a smooth, non-degenerate $k$-normal 
surface, whose hyperplane section is denoted by $H$. Assume that $h^1 (\os(kH))=0$. 

Let $|D|$ be a Bertini linear system on $S$ whose general element is $k$-normal.
Assume that the following inequalities hold:
$$D.H>kH^{2}=kdeg(S), \quad (D-2kH)^{2}>0,\quad D(D-2kH)>0,$$ 
$$\nu(D,kH)=k^{2}((D.H)^{2}-D^{2}H^{2})<4(D.(D-2kH)-4),$$
where $\nu(D,kH)$ is the Hodge number of $D$ and $kH$, and finally:
$$\delta<f(H,D,S):=\frac{C.(C-2kH)+\sqrt{C^{2}(C-2kH)^{2}}}{8}.$$
Thus, if $C\in |D|$ is a reduced irreducible curve having $\delta$ nodes and cusps as 
the only singularities, then $C$ is geometrically $k$-normal.
}

\vskip 5pt

In particular, we have:

\vskip 5pt
\noindent
{\bf Corollary} (cfr. Corollary \ref{exampleprop}) 
{\em Let $S$ be a smooth non-degenerate projective complete intersection surface 
in $\mathbb{P}^{r}$, $r\geq 3$, and let
$C\in |nH|$ be an irreducible curve having only $\delta$ nodes and cusps as singular points. 
Assume $n\geq 2k+1$, $deg(S)>\frac{4}{n(n-2k)}$; if 
$$\delta<\frac{n(n-2k)deg(S)}{4}$$then $C$ is g.k.n..
}

\vskip 5pt
\noindent
The strategies for proving the above results are 
similar to those used in \cite{fm} for proving analogous theorems.

Furthermore, in Section 2 we investigate on some 
possible relations between geometric $k$-normality and geometric
$(k-1)$-normality of curves (cfr. Proposition \ref{gkng(k-1)n}).

In Section 3 we relate the notion of geometric $k$-normality to other 
important properties of a singular curve $C$ on a smooth surface $S$ 
and of its normalization $\tilde{C}$. 
A first instance of this connection has already been studied in \cite{flaminio},
\S 5, in the case of geometric linear normality. Here we make some further 
observations for the case of geometric $k$-normality. Indeed, if $N= Sing(C)$, 
we show that the g.k.n. property of $C$ is related, on the one hand, to the 
$0$-regularity (in the sense of Castelnuovo-Mumford)
of a suitable twist of the ideal sheaf $\mathcal{I}_{N}$ (cfr. Lemma \ref{olemma}) 
on the other hand, to other intrinsic properties of $\tilde{C}$ such as its {\em Brill-Noether 
number} (cfr. Proposition \ref{mu2}).

In particular, we have:

\vskip 5pt

\noindent
{\bf Proposition} (cfr. Corollary \ref{mu3}) {\em Assume that 
$S$ is a smooth linearly normal and 2-normal, non-degenerate 
surface in $\mathbb{P}^{r}$, $D-3H$ is big and nef and that 
$h^1(\os)=h^1(\os(H))=h^1(\os(2H))=0$. Let $C\in |D|$ be a reduced, irreducible 
curve having only $\delta$ nodes and cusps as singularities and let 
$\phi:  \tilde{C} \to C$ be its normalization. Let $g$ be the geometric 
genus of $C$ (and so of $\tilde{C} $).

If $C$ is geometrically 2-normal, 
then the Brill-Noether map $\mu_{0, \tilde{C}}$ of the pair $(\tilde{C}, \phi^{*}H)$ 
is surjective; in particular the Brill-Noether number 
$$\rho(g,r,deg(C)):=g-(r+1)(r-deg(C)+g)\leq 0.$$Thus, if $\rho(g,r,deg(C))>0$, under the 
hypotheses above, $C$ cannot be geometrically 2-normal. 
}

\vskip 5pt

\noindent
We also discuss some interesting examples of geometrically linearly normal 
but not geometrically 2-normal curves on smooth quadrics  in $\mathbb{P}^3$ and 
which have positive Brill-Noether number 
(cfr. Example \ref{ex:vai}). 

In Section 4 we show how the concept of geometric $k$-normality 
is also closely related to the study of 
some {\em Severi varieties} which parametrize families of 
irreducible and $\delta$-nodal curves on smooth surfaces. 

We first focus on surfaces of general type and 
we show that the bound on $\delta$ in the Theorem above is sharp (at least 
for curves $C\in |8H|$ of degree 48, with 48 nodes and living on a 
smooth sextic surface in $\mathbb{P}^3$). From the relation between geometric
$k$-normality and the local behaviour of Severi varieties, it turns out that 
in the above example the point $[C]\in \mathcal{V}_{48,|8H|}$ is not a {\em regular 
point} of the given Severi variety (for precise definitions cf. \S 4). 

Finally, we consider the case 
of Severi varieties of 
$\delta$-nodal curves of degree $n$ in the plane and we give an 
upper-bound on the number of nodes $\delta\leq f(n,k)$ (valid for $k=1, \; 2,\;3$)
such that if $\delta\leq f(n,k)$, then the general element of the 
corresponding Severi variety $V_{n, \delta}$ represents a geometrically $k$-normal plane 
curve (cfr . Theorem \ref{sv}). In the case $k=1$, the result was proved 
by Sernesi in \cite{ser}. 
    
We use standard notation throughout the paper. If $\mathcal{F}$ is a coherent 
sheaf on some projective variety $X$, $H^i(X,\mathcal{F})$ will denote the i-th
cohomology group and $h^i(X,\mathcal{F})$ the corresponding dimension. If
$Z \subset X$ is any divisor in $X$, $|Z|$ will denote the complete linear 
system on $X$ corresponding to $Z$, i.e. the projective space parametrizing 
effective divisors linearly equivalent to $Z$. The symbol $Z \sim Z'$ will 
denote linear equivalence. If $C$ is a reduced irreducible (possibly singular) 
curve, $p_a(C)$ will denote its arithmetic genus, whereas $g(\tilde{C})$ will 
denote its geometric genus, i.e. the arithmetic genus of its normalization
$\tilde{C}$.  

\section{Geometric $k$-normality on some projective surfaces}
In this section we study under which conditions an 
irreducible nodal and cuspidal curve $C$ on a projective surface $S$ is 
geometrically $k$-normal.
The strategies for proving Theorem \ref{thm1} and Theorem \ref{thm2} are
similar to those used in \cite{cs} and \cite{fm} for analogous results; 
nevertheless we give detailed proofs of Thereom \ref{thm1} and Theorem 
\ref{thm2} not only for convenience of a reader, but mainly because we 
consider the more general case of curves with cusps and nodes as 
singularities. 

The first result is the following, which is a direct generalization of Theorem
1 in \cite{fm}.
\begin{theorem}\label{thm1}
Let $S$ be a smooth non-degenerate $k$-normal surface in $\mathbb{P}^r$, 
such that $h^1(S,\mathcal{O}_S(kH))=0$ (e.g. $S$ a complete intersection), with canonical 
divisor class $K_S$ and hyperplane divisor class $H$. 
Let $|D|$ be a Bertini linear system on $S$ such that its general element is 
$k$-normal (e.g. $D\sim nH$ on $S$ complete intersection). 
Let $C \in |D|$ be an irreducible curve having only nodes and cusps as singularities, let $N=Sing(C)$ be 
the reduced zero dimensional subscheme in $S$ of the singularities of $C$ and let $\delta$ be the length of $N$; 
then $C$ is g.k.n. iff $N$ imposes indipendent conditions on the linear system $|D+K_S-kH|$, i.e.
$h^0(\mathcal{I}_{N}(D+K_S-kH))=h^0(\mathcal{O}_S(D+K_S-kH))-\delta$.
\end{theorem}
\begin{proof} 
First of all, let us remark that under the hypotheses of Thereom \ref{thm1} $C$ is g.k.n. iff, $h^0(\mathcal{O}_{\tilde{C}}(\phi^{*}(kH)))=h^0(\mathcal{O}_D(kH))=\binom{r+k}{r}-P(k)$, where $P(k)$ is the postulation in degree $k$ of $D$ 
and $\phi: \tilde{C} \rightarrow C$ is the normalization map. 
This is an immediate consequence of Definition \ref{gkn}, combined with the assumption that $D$ is $k$-normal.
Let $D$ be the general element of the linear system $|D|$. By using 
Riemann-Roch and the $k$-normality hypothesis we have 
$h^1(\mathcal{O}_D(kH))=\binom{r+k}{r}-kdeg(D)+p_a(D)-1-P(k)$; 
by Serre's duality and by the adjunction formula on $S$, we get:
\begin{equation}\label{one}
h^0(\mathcal{O}_D(D+K_S-kH))=h^0(\mathcal{O}_D(kH))-kdeg(D)+p_a(D)-1.
\end{equation}
Let $\mu: \tilde{S}\rightarrow S$ be the blowing-up of $S$ along the reduced 
zero dimensional subscheme of singularities $N$ and let $B=\sum_{i=1}^{\delta}E_i$ be the exceptional divisor. The restriction of $\mu$ to the strict transform $\tilde{C}$ of $C$ in $\tilde{S}$ is the normalization map $\phi$ (since we are admitting only nodes or cusps as singularities).
By adjunction on $\tilde{S}$, we have 
\begin{equation}\label{two}
\omega_{\tilde{C}}=
\mathcal{O}_{\tilde{C}}(K_{\tilde{S}}+\tilde{C})=\ot(\mu^{*}(K_S+C)-B)
=\ot(\phi^{*}(K_S+C)-\tilde{N}),
\end{equation}
where $\tilde{N}$ is an 
effective divisor of degree $2\delta$ on $\tilde{C}$ which maps to the singularities of $C$.
By Riemann-Roch on $\tilde{C}$, we get: $h^1(\tilde{C},\ot(\phi^{*}(kH)))=
h^0(\tilde{C},\ot(\phi^{*}(kH)))-kdeg(C)+g(\tilde{C})-1$, where $g(\tilde{C})=
p_a(C)-\delta$.
By Serre's duality and by using \eqref{two} and the fact that $C\sim D$ on $S$, we get:
$h^1(\ot(\phi^{*}(kH)))=h^0(\ot(\phi^{*}(K_S+D-kH))-\tilde{N})$. In this way we get:
\begin{equation}\label{three}
h^o(\ot(\phi^{*}(K_S+D-kH)-\tilde{N}))=h^0(\ot(\phi^{*}(kH)))-kdeg(C)+p_a(C)-1-\delta.
\end{equation}
Now, observe that $h^0(\tilde{C},\ot(\phi^{*}(kH)))=\binom{r+k}{r}-P(k)=h^0(D,
\mathcal{O}_D(kH))$, i.e. $C$ is g.k.n. iff, by \eqref{three},
\begin{equation}\label{four}
h^o(\ot(\phi^{*}(K_S+D-kH)-\tilde{N}))=h^0(\mathcal{O}_{D}(kH))-kdeg(C)+p_a(C)-1-\delta.
\end{equation}
Combining \eqref{one}, \eqref{four} and the fact that adjunction on $S$ is 
independent of the chosen element in $|D|$, we get that $h^0(\ot(\phi^{*}(kH)))=h^0(\mathcal{O}_{D}(kH))$ (i.e. $C$ is g.k.n.) iff:
\begin{equation}\label{five}
h^0(\ot(\phi^{*}(K_S+D-kH)-\tilde{N})=h^0(\mathcal{O}_{D}(K_S+D-kH))-\delta=
\end{equation}
\[  =h^0(\mathcal{O}_{C}(K_S+C-kH))-\delta.\]
Thus, using the assumption $h^1(S,\mathcal{O}_{S}(kH))=0$, Serre's duality on $S$ and the exact sequence
\[
0 \rightarrow \mathcal{O}_{S}(K_S-kH) \rightarrow \mathcal{O}_{S}(C+K_S-kH) 
\rightarrow \mathcal{O}_{C}(C+K_S-kH) \rightarrow 0
\]
we obtain:
\begin{equation}\label{six}
h^0(\mathcal{O}_{C}(C+K_S-kH))=h^0(\mathcal{O}_{S}(K_S-kH+C))-h^0(\mathcal{O}_{S}(K_S-kH)).
\end{equation}
Analogously, from the exact sequence on $\tilde{S}$:
\[
0 \rightarrow \ost(K_{\tilde{S}}-\mu^{*}(kH)) \rightarrow 
\ost(\tilde{C}+K_{\tilde{S}}-\mu^{*}(kH)) \rightarrow \ot(\tilde{C}+K_{\tilde{S}}-\mu^{*}(kH))\rightarrow 0
\]
we get the following exact sequence:
\[
0 \longrightarrow \ost(\mu^{*}(K_S-kH)+B) \longrightarrow \ost(\mu^{*}(C+K_S-kH)-B) \longrightarrow\]
\[ \longrightarrow \ot(\phi^{*}(C+K_S-kH)-\tilde{N}) \longrightarrow 0
\]
from which, we obtain
\begin{equation}\label{seven}
h^0(\ot(\phi^{*}(C+K_S-kH)-\tilde{N}))=
\end{equation}
\[=h^0(\ost(\mu^{*}(C+K_S-kH)-B))
-h^0(\ost(\mu^{*}(K_S-kH)+B)). 
\]
Indeed, applying Serre's duality on $\tilde{S}$, 
we get the equality $h^1(\ost(\mu^{*}(K_S-kH)+B))=h^1(\ost(\mu^{*}(kH)))$ and this is zero since $h^1(\ost(\mu^{*}(kH)))=h^1(\os(kH))=0$ as it follows from 
the (degenerate) Leray spectral sequence.
Combining \eqref{five}, \eqref{six} and \eqref{seven} we find:
\begin{equation}\label{eight}
h^0(\ost(\mu^{*}(C+K_S-kH)-B))-h^0(\ost(\mu^{*}(K_S-kH)+B))=
\end{equation}
\[=h^0(\os(K_S-kH+C))-h^0(\os(K_S-kH))-\delta.
\]
Now observe that applying Serre's duality on $\tilde{S}$ and on $S$, and using the (degenerate) Leray spectral sequence we get that $h^0(\os(\mu^{*}(K_S-kH)+B))=
h^2(\os(\mu^{*}(kH)))=h^2(\os(kH))=h^0(\os(K_S-kH))$. 
Thus we end up with the desired result using \eqref{eight} and observing that 
$h^0(\ost(\mu^{*}(C+K_S-kH)-B)=h^0(\mathcal{I}_{N}(D+K_S-kH))$.
\end{proof}

\begin{remark}
In the previous theorem we have implicitly assumed that $|D+K_S-kH|\neq
\emptyset$. If this is not the case, by Riemann-Roch we have that $C$ is g.k.n.
 iff $\delta=0$, i.e. $C$ is smooth.   
\end{remark}

One of our goals is to apply Theorem \ref{thm1} to the study of 
Severi varieties in 
\S \ref{severivariety}, so we have restricted our consideration 
just to the case of irreducible nodal and cuspidal curves. 
However, we strongly believe that Theorem \ref{thm1} admits a direct 
generalization to irreducible curves having {\em any type} of singularities,
by using the adjoint/conductor ideal sheaf of $C$. This generalization will be 
not be considered here and it will be the subject of a future analysis.

Anyway, the following proposition is immediate:
\begin{proposition}
Let $S$ and $|D|$ be as in Theorem \ref{thm1}. Let $C\in |D|$ be an irreducible curve 
having only $\alpha$ ordinary singular points of multiplicity $m$. 
Let $\Sigma \subset |D+K_S-kH|$ be the linear system of 
curves having singular points of multiplicity $m-1$ at the 
singular points of $C$. Then $C$ is g.k.n. iff $dim(\Sigma)=dim(|D+K_S-kH)|)-
\alpha \frac{m(m-1)}{2}$. 
\end{proposition}
\begin{proof}
It is enough to observe that the blow-up of $S$ at the singular points of $C$,
$\mu: \tilde{S} \rightarrow S$ restricts to the normalization map $\phi: \tilde{C} \rightarrow C$ on the strict transform $\tilde{C}$ of $C$ in $\tilde{S}$.
Then the same proof given for Theorem \ref{thm1} will do the job here.
\end{proof}
Next, we find sufficient and numerical conditions on $D$, $S$ and $H$ which 
determine upper-bounds $f(H,D,S)$ on $\delta$, such that,
if $C\in |D|$ is a curve 
as in Theorem \ref{thm1} and $\delta<f(H,D,S)$, 
then the subscheme of singularities $N$ imposes independent conditions on 
$|D+K_S-kH|$. Following again \cite{fm} and \cite{cs}, the idea is to construct a Bogomolov-unstable 
vector bundle on $S$ to obtain the function $f(H,D,S)$, so that if the number of nodes and cusps $\delta$ is 
at most $f(H,D,S)-1$, then $C$ is geometrically $k$-normal.
In the proof of this result we shall make use of the following:
\begin{definition}\label{bu}
 Let $S$ be a smooth projective surface; a rank 2 vector bundle $\mathcal{E}$ on $S$ is called {\em Bogomolov-unstable} if there exist $M,B 
\in Div(S)$ and a 0-dimensional subscheme $Z$ (possibly empty) such that 
$\mathcal{E}$ fits in the following exact sequence:
\[
0 \rightarrow \os(M) \rightarrow \mathcal{E} \rightarrow \mathcal{I}_{Z}(B) 
\rightarrow 0 
\]
and moreover $(M-B)\in N(S)^{+}$, where $N(S)^{+}$ is the ample cone on $S$.
Let us also recall that if $c_1(\mathcal{E})^{2}-4c_2(\mathcal{E})>0$, then $\mathcal{E}$ is Bogomolov-unstable (see for instance \cite{friedman}). 
\end{definition}
We can now state the following:
\begin{theorem}\label{thm2}
Let $S\subset \mathbb{P}^{r}$ be a smooth non-degenerate surface, whose hyperplane 
section is denoted by $H$. Let $D$ be an
effective irreducible smooth divisor on $S$ and let $k$ be a positive integer. 
Assume that the following 
inequalities hold:
\begin{equation}\label{2.1}
D.H>kH^{2}=kdeg(S),
\end{equation}
\begin{equation}\label{2.2}
(D-2kH)^{2}>0,\quad D(D-2kH)>0, 
\end{equation}
\begin{equation}\label{2.3}   
\nu(D,kH)=k^{2}((D.H)^{2}-D^{2}H^{2})<4(D.(D-2kH)-4),
\end{equation}
where $\nu(D,kH)$ is the Hodge number of $D$ and $kH$, and finally:
\begin{equation}\label{2.4}
\delta<f(H,D,S):=\frac{C.(C-2kH)+\sqrt{C^{2}(C-2kH)^{2}}}{8}.
\end{equation}
Then if $C\in |D|$ is a reduced irreducible curve having only $\delta$ nodes and cusps as singularities, we have that $N=Sing(C)$ imposes independent conditions on $|D+K_S-kH|$.
\end{theorem}
\begin{proof}
Assume the contrary, that is $N$ does not impose independent conditions on 
$|D+K_S-kH|$. Let $N_{0}\subset N$ be the 0-dimensional subscheme of $N$ of minimal length for which $N_{0}$ does not impose independent conditions on 
$|D+K_S-kH|$ and let $\delta_{0}$ be the length of $N_{0}$. Due to the exact 
sequence:
\[ 0 \rightarrow \mathcal{I}_{N_{0}}\otimes\os(D+K_S-kH) \rightarrow
\os(D+K_S-kH) \rightarrow \underline{\mathbb{C}}^{N_{0}} \rightarrow 0,
\]
the failure of $N_{0}$ to impose independent coditions means $h^1(\mathcal{I}_{N_{0}}\otimes \os(D+K_S-kH))\neq 0$ and so $N_{0}$ satisfies the Cayley-Bacharach condition (i.e. if an element of the linear system $|D+K_S-kH|$ passes through all the components of $N_{0}$ except one, then is has to pass also through the remaining one).
Therefore (see for instance \cite{la}) a non-zero element of $H^{1}(
\mathcal{I}_{N_{0}}(D+K_S-kH))$ gives rise to a non-trivial rank 2 vector 
bundle $\mathcal{E}\in \ext^{1}(\mathcal{I}_{N_{0}}(D-kH), \os)$, fitting in 
the following exact sequence:
\begin{equation}\label{2.5}
0 \rightarrow \os \rightarrow \mathcal{E} \rightarrow \mathcal{I}_{N_{0}}(D-kH) \rightarrow 0,
\end{equation}  
so that $c_1(\mathcal{E})=C-kH$ and $c_2(\mathcal{E})=\delta_{0}$.
In our case, we have:
\begin{equation}\label{2.6}
c_1(\mathcal{E})^{2}-4c_2(\mathcal{E})=(D-kH)^{2}-4\delta_{0}.
\end{equation}
Now by \eqref{2.4} we have $(D-kH)^{2}-4\delta_0\geq (D-kH)^{2}-4\delta$ which is   grater than $k^{2}H^{2}+\frac{1}{2}D.(D-2kH)>k^{2}H^{2}>0$, 
since, by \eqref{2.2}, $D.(D-2kH)>0$.
Thus, $\mathcal{E}$ is Bogomolov-unstable and hence, using Definition \ref{bu}, $h^0(\mathcal{E}(-M))\neq 0$. Twisting \eqref{2.5} by $\os(-M)$ we get:
\begin{equation}\label{2.7}
0 \rightarrow \os(-M) \rightarrow \mathcal{E}(-M) \rightarrow \mathcal{I}_{N_{0}}(D-kH-M) \rightarrow 0.
\end{equation}
Observe that $h^{0}(\os(-M))=0$: otherwise, $-M$ would be an effective divisor and we would have $-M.A>0$ for each ample divisor $A$. Then by the exact sequence in Definition \ref{bu} it follows that $c_1(\mathcal{E})=M+B$, and so, by the exact sequence \eqref{2.5}, 
\begin{equation}\label{2.8}
M-B=2M-D+kH \in N(S)^{+}.
\end{equation}
Thus $(2M-D+kH).H>0$, which implies
\begin{equation}\label{2.9}
M.H>\frac{(D-kH).H}{2},
\end{equation}
so that, by \eqref{2.1}, $H.(D-kH)>0$, and we get $-M.H<0$, which gives the claim $h^{0}(\os(-M))=0$.
From the long exact cohomology sequence of \eqref{2.7}, we get $0\neq H^0(
\mathcal{E}(-M)) \hookrightarrow H^0(\mathcal{I}_{N_{0}}(D-kH-M))$. This means that there exists $\Delta \in |D-kH-M|$ such that $N_0\subset \Delta$. If it 
were $C\subset \Delta$, then $\Delta-C\geq 0$, but on $S$ we have 
$\Delta-C \sim -kH-M$, so that also $-kH-M$ would be effective, but this is clearly false since $-kH^2-H.M<0$ by \eqref{2.9}. 
Thus $C$ cannot be a component of $\Delta$. Applying Bezout's theorem we get:
\begin{equation}\label{2.10}
C.\Delta=C.(D-kH-M)\geq 2\delta_{0},
\end{equation}
since $C$ has only nodes and cusps as singularities.
On the other hand, taking $M$ maximal among all divisors satisfying $h^0(
\mathcal{E}(-M))\neq 0$, we can further assume that the general section of 
$\mathcal{E}(-M)$ vanishes in codimension 2. Denoting by $Z$ this 
vanishing locus, we have $c_2(\mathcal{E}(-M))=deg(Z)\geq 0$. But $c_2(
\mathcal{E}(-M))=c_2(\mathcal{E})+M^{2}+c_1(\mathcal{E}).(-M)=\delta_{0}+M^{2}-
M.(D-kH)$, which implies
\begin{equation}\label{2.11}
\delta_{0}\geq M.(D-kH-M).
\end{equation}
Applying the Hodge Index Theorem to the divisor pair $(D,2M-D+kH)$ we 
obtain:
\begin{equation}\label{2.12}
D^{2}(2M-D+kH)^{2}\leq(D.(2M-D+kH))^{2}=
\end{equation}
\[
=(D.(D-kH)-2D.(D-kH-M))^{2}.
\]
From \eqref{2.1}, and the second part of \eqref{2.2}  it follows $D.(D-kH)>0$ 
and similarly one finds $D^{2}>0$ (this is also due to our assumption that $|D|$ is a mobile linear system). Since $D$ is irreducible, this implies that $D$ is nef.
From \eqref{2.10} and the positivity of $D.(D-kH)$, it follows that
\begin{equation}\label{2.13}
D.(D-kH)-2D.(D-kH-M)\leq D.(D-kH)-4 \delta_{0}.
\end{equation}
Next observe that the left hand side of \eqref{2.13} is non-negative, since it is equal to $D.(2M+D-kH)$, where $D$ is effective and, by \eqref{2.8}, we have 
$2M-D+kH \in N(S)^{+}$. Squaring both sides of \eqref{2.13} and comparing with 
\eqref{2.12} we get
\begin{equation}\label{2.14}
D^{2}(2M-D+kH)^{2}\leq (D.(D-kH)-4 \delta_{0})^{2}.
\end{equation}
On the other hand, we find: $(2M-D+kH)^{2}=
4(M-\frac{D-kH}{2})^{2}$, which is equal to $(D-kH)^{2}-4(D-kH+M).M$ and, using \eqref{2.11},
\begin{equation}\label{2.15}
(2M-D+kH)^{2}\geq (D-kH)^{2}-4 \delta_{0}.
\end{equation}
Next define the following function 
\begin{equation}\label{2.16}
F(\delta_{0}):=16\delta_{0}^{2}-4\delta_{0}D.(D-2kH)
+k^{2}((D.H)^{2}-D^{2}H^{2}).
\end{equation}
Using \eqref{2.14} and \eqref{2.15}, it turns out that $F(\delta_{0})\geq 0$. 
Now we will show that under our numerical hypotheses, one has 
$F(\delta_{0})<0$ for some positive integer $\delta_{0}$, so that the 
assumption $h^1(\mathcal{I}_{N}(D-kH+K_S))\neq 0$ leads to a contradiction. 
To do this, consider the reduced discriminant of the equation $F(\delta_{0})=0$; 
this is equal to $4D^2(D-2kH)^2$, which is strictly positive due to the first
part of \eqref{2.2} and to the fact that $D^2>0$. Thus $F(\delta_{0})<0$ iff $\delta_0 \in 
(\alpha(D,kH), \beta(D,kH))$ where $\alpha(D,kH)$ and $\beta(D,kH)$ are the 
two roots of the equation. It is a straightforward computation to see 
that $\delta_{0}<\beta(D,kH)$ is just condition \eqref{2.4}. On the other hand 
$\alpha(D,kH)\geq 0$: if $\alpha(D,kH)<0$, then it turns out that $D.(D-2kH)<\sqrt{D^2(D-2kH)^2}$, which contradicts 
the Hodge Index Theorem since $D.(D-2kH)>0$ by 
\eqref{2.2}.

Finally, we have to show that there exists at least one positive integer value 
of $\delta_{0}$ between $\alpha$ and $\beta$. In order to simplify the 
notation, put $t:=D.(D-2kH)$. Thus, from the expression of $\alpha(D,kH)$ it 
turns out that $\alpha<1$ iff $t-8<\sqrt{t^2-4\nu(D,kH)}$. If $t-8<0$, the previous inequality trivially holds, so that $\delta_{0}>\alpha(D,kH)$. In this case note also that by \eqref{2.3} we find $\beta(D,kH)>1$, so that there exists at least one positive integral value for the number of nodes and cusps 
$\delta_{0}$.
On the other hand, if $t-8 \geq 0$, then $\alpha(D,kH)<1$ directly from 
\eqref{2.3}, while $\beta(D,kH)>1$ holds since it is equivalent to $t-8>-
\sqrt{t^2-4\nu(D,kH)}$. 
\end{proof}  
Combining Theorems \ref{thm1} and \ref{thm2} we immediately get the following:
\begin{corollary}\label{aiuto}
Let $S$ be a smooth non-degenerate $k$-normal surface in $\mathbb{P}^r$, 
such that $h^1(\os(kH))=0$. 
Let $|D|$ be a Bertini linear system on $S$ whose general element is $k$-normal.  Assume that $D$ satisfies the 
inequalities \eqref{2.1}, 
\eqref{2.2} and \eqref{2.3} of Theorem \ref{thm2}. If $C \in |D|$ is a reduced irreducible curve having only $\delta$ 
nodes and cusps as singularities and $\delta$ satisfies \eqref{2.4}, then $C$ is geometrically $k$-normal. 
\end{corollary}

\begin{remark}
As it is clear from the proof, in Theorem \ref{thm2} we just find some 
numerical conditions ensuring that $h^1(\mathcal{I}_{N}(D+K_S-kH))=0$. Observe, however, that 
in general this is only a {\em sufficient} condition for the subscheme $N$ of nodes and cusps to impose 
independent conditions on the linear system $|D+K_S-kH|$. However, by the Kawamata-Viehweg Vanishing Theorem, 
the condition $h^1(\mathcal{I}_{N}(D+K_S-kH)) = 0$ is also necessary if $D-kH$ is nef and big. 
\end{remark}

In the next sections, we will consider some examples of geometrically $k$-normal nodal curves 
on smooth complete intersection surfaces. 
In order to study these examples, we first rewrite the conditions of 
Theorems \ref{thm1} and \ref{thm2} for surfaces which are 
smooth complete intersections in a 
more explicit form. 

We have the following: 
\begin{corollary}\label{exampleprop}
Let $S$ be a smooth non-degenerate projective surface, which is a 
complete intersection in $\mathbb{P}^{r}$, $r\geq 3$, and let
$C\in |nH|$ be an irreducible curve having only $\delta$ nodes and cusps as singular points. 
Assume $n\geq 2k+1$, $deg(S)>\frac{4}{n(n-2k)}$; if 
\begin{equation}\label{exampleprop1}
\delta<\frac{n(n-2k)deg(S)}{4}
\end{equation}
then $C$ is g.k.n..
\end{corollary}

\begin{proof}
The assumption $n\geq 2k+1$ is just \eqref{2.1} and \eqref{2.2} rewritten for this case. 
The bound on the degree of $S$ is obtained from \eqref{2.3}, whereas the bound on the 
number of nodes and cusps \eqref{exampleprop1} is obtained from 
\eqref{2.4}, recalling that if the Hodge number is zero as in this case is, 
then $\sqrt{C^2(C-2kH)^2}=C(C-2kH)$. Finally, to apply Theorem \ref{thm1} and 
to conclude that $C$ is geometrically $k$-normal we must prove that the general element 
of $|nH|$ is k-normal. This is immediate.
\end{proof}

To conclude this section, we investigate on possible relations between geometric 
$k$-normality and geometric $(k-1)$-normality. A first insight is given by the following:
\begin{proposition}\label{gkng(k-1)n} 
Let $S$, $|D|$ and $C\in|D|$ as in the hypotheses of Theorem \ref{thm1}; 
assume moreover that $h^1(\os(kH))=h^1(\os((k-1)H))=0$ and that $S$ is $k$-normal and $(k-1)$-normal. Then if $C$ 
is geometrically 
$k$-normal, it is also geometrically $(k-1)$-normal. 
\end{proposition}

\begin{proof}
By Theorem \ref{thm1}, $C\in|D|$ is g.k.n. iff the subscheme $N$ of singularities of $C$ (consisting of 
only nodes and cusps) imposes independent conditions on $|D+K_S-kH|$. This is clearly equivalent to claim that 
the evaluation map $\psi_k : H^0(\os(D+K_S-kH)) \rightarrow \underline{\mathbb{C}}^{\delta}$ is surjective, where 
$\delta$ is the length of $N$ and $\underline{\mathbb{C}}^{\delta}$ is a 
skyscraper sheaf supported on $N$. On the other hand, $\psi_k$ is surjective 
iff there exist sections $\sigma_i \in H^0(\os(D+K_S-kH))$ such that 
$\sigma_i(P_j)=\delta_{i,j}$ (where $\delta_{i,j}$ is the Kronecker's symbol), 
$i,j=1\dots\delta$, where $P_1,\dots,P_{\delta}$ are the singularities of $C$. Choose 
a section $\sigma \in H^0(S, \os(H))$ such that $\sigma(P_j)\neq 0$ for any $j$. Now consider the 
multiplication map $$m :H^0(\os(H))\otimes H^0(\os(D+K_S-kH))\rightarrow
H^0(\os(D+K_S-(k-1)H)).$$ The sections $g_i:=m(\sigma\otimes\sigma_i)$, $i=1,\dots, \delta$, belong to 
$H^0(\os(D+K_S-(k-1)H))$ and they satisfy $g_i(P_j)=
\delta_{i,j}$. Using the $g_i$'s it is immediate to see that if $\psi_{k-1}:
H^0(\os(D+K_S-(k-1)H)) \rightarrow \underline{\mathbb{C}}^{\delta}$ is the valuation map related to $|D+K_S-(k-1)H|$, 
then $\psi_k$ surjective implies $\psi_{k-1}$ surjective.    
\end{proof}

In the next section, we 
observe that the converse of the above result 
does not hold (cf. Example \ref{ex:vai})

\section{Intrinsic and extrinsic properties of geometric $k$-normality}

The notion of geometric $k$-normality is strictly connected to other important
topics, such as the 0-regularity of Castelnuovo-Mumford of some 
suitable coherent sheaves on $S$, and also to geometric intrinsic properties 
of smooth curves, e.g. the well-known 
{\em Brill-Noether number}. A first instance of this 
connection has already been studied in \cite{flaminio}, \S 5. Here we make 
some further observations on this subject.  

We start to analyze the relation between geometric k-normality of $C\subset S$ 
and the 0-regularity (in the sense of Castelnuovo-Mumford) of the coherent sheaf 
$\mathcal{I}_{N}(D+K_S-(k-1)H)$. Recall that a coherent sheaf $\mathcal{F}$ 
on $S$ is said to be {\em $k$-regular} iff $h^i(S,\mathcal{F}(k-i))=0$ for any $i>0$.

This relation is motivated by the following 
observation: under the hypotheses of Theorem \ref{thm1} we know that 
$h^1(\mathcal{I}_{N}(D+K_S-kH))=0$ implies that $C$ is geometrically 
$k$-normal. Now, if $D-kH$ is ample (or big and nef), then by the 
Kawamata-Viehweg vanishing theorem we have $h^1(\os(K_S+D-kH))=0$. In this case, we have that $C$ is 
g.$k$.n. iff $h^1(\mathcal{I}_{N}(D+K_S-kH))=0$; (for instance, this is the case if $S$ is a c.i. 
and $D\sim nH$, with $n>k$). For what concerns the 
0-regularity of $\mathcal{I}_{N}(D+K_S-kH)$, we observe the following:

\begin{lemma}\label{olemma}
If $D-(k+1)H$ is big and nef, then $h^1(\mathcal{I}_{N}(D+K_S-kH))=0$ iff
$\mathcal{I}_{N}(D+K_S-(k-1)H)$ is 0-regular. Moreover, if we assume $h^1(\os(kH))=0$, then both 
conditions are equivalent to saying that $C\in |D|$, with 
only $\delta$ nodes and cusps as singularities, is geometrically $k$-normal.
\end{lemma} 

\begin{proof}
By the exact sequence
\[0 \rightarrow \mathcal{I}_{N}(K_S+D-(k+1)H) \rightarrow \os(K_S+D-(k+1)H)
\rightarrow \mathcal{O}_N \rightarrow 0 
\]
and using Kawamata-Viehweg vanishing we immediately get 
$h^2(\mathcal{I}_{N}(D+K_S-(k+1)H))=0$. The rest comes from 
Theorem \ref{thm1} and from what observed above.
\end{proof}

Let us now recall the following well-known proposition (for a proof see \cite{m}):
\begin{proposition}\label{castelnuovo}
Let $\mathcal{F}$ be a coherent sheaf on a smooth projective variety $X$. Then if $\mathcal{F}$ is m-regular on $X$, we have that 
\begin{itemize}
\item 
the maps:
\[ H^0(X,\mathcal{F}(l-1))\otimes H^0(X, \mathcal{O}_X(1)) \rightarrow H^0(X, \mathcal{F}(l))
\]
are surjective for any $l>m$; 
\item $H^i(X,\mathcal{F}(l))=0$, for $i>0$, $l+i\geq m$.
\end{itemize}
\end{proposition}

Observe that the 0-regularity of the coherent sheaf 
$\mathcal{I}_{N}(K_S+D-(k-1)H)$ has some interesting consequences on the behaviour 
of linear systems on $\tilde{C}$. 
Indeed the 0-regularity implies that
\[
H^0(\mathcal{I}_{N}(K_S+D-(k-1)H))\otimes H^0(\os(H))\rightarrow
H^0(\mathcal{I}_{N}(D+K_S-(k-2)H))
\]
is surjective, by the first part of Proposition \ref{castelnuovo}.
This is equivalent to the surjectivity of the following map:
\[
H^0(\tilde{S}, \ost(\mu^{*}(D+K_S-(k-1)H)-B)) \otimes H^0(\tilde{S}, \ost(
\mu^{*}(H))) \rightarrow 
\]
\[ 
\rightarrow  H^0(\ost(\mu^{*}(K_S+D-(k-2)H)-B)).
\]Moreover, if we assume $h^1(\os((k-1)H))=h^1(\os((k-2)H))=0$, 
then, reasoning  as in \cite{flaminio} (page 747), we find that 
\begin{equation}\label{mu1}
\mu_{k-2, \tilde{C}}: H^0(\omega_{\tilde{C}}(-(k-1)\mu^{*}(H))) \otimes
H^0(\ot(\mu^{*}(H))) \rightarrow H^0(\omega_{\tilde{C}}(-(k-2)H)),
\end{equation}
is also surjective.
Thus we have proved the following 
\begin{proposition}\label{mu2}
The 0-regularity of the sheaf $\mathcal{I}_{N}(K_S+D-(k-1)H)$ and 
the vanishing of the cohomology groups $H^1(\os((k-1)H))$ and 
$H^1(\os((k-2)H)$ imply the surjectivity of \eqref{mu1}. In particular, for 
$k=2$, this gives the surjectivity the well-known Brill-Noether map
$$\mu_{0, \tilde{C}}:H^0(\omega_{\tilde{C}}(-\phi^{*}(H)))\otimes 
H^0(\mathcal{O}_{\tilde{C}}(\phi^{*}(H))) \rightarrow H^0(\omega_{\tilde{C}}).$$
\end{proposition}

\begin{corollary}\label{mu3}
Assume that $S$ is a smooth linearly normal and 2-normal non-degenerate 
surface in $\mathbb{P}^{r}$, $D-3H$ is big and nef and that 
$h^1(\os)=h^1(\os(H))=h^1(\os(2H))=0$. If $C\in |D|$ is geometrically 2-normal, 
then the Brill-Noether map $\mu_{0, \tilde{C}}$ of the pair $(\tilde{C}, \phi^{*}H)$ 
is surjective; in particular the Brill-Noether number 
$$\rho(g,r,deg(C)):=g-(r+1)(r-deg(C)+g)\leq 0.$$Thus, if $\rho(g,r,deg(C))>0$, under the 
hypotheses above, $C$ cannot be 
geometrically 2-normal. 
\end{corollary} 
\begin{proof}
Observe that under the above hypotheses on $S$, the geometric 2-normality of $C$ 
implies its linear normality by Proposition \ref{gkng(k-1)n}. Thus if $C$ is also 
geometrically linearly normal, then the Brill-Noether number can be 
expressed as the difference between the dimension of the cokernel and the 
kernel of the Brill-Noether map. 
\end{proof}

\begin{example}\label{ex:vai}
Let us discuss the following example in order to clarify the meaning of Proposition \ref{mu2} and Corollary  \ref{mu3} and their limits. 
Consider a smooth quadric $S_{2}\subset \mathbb{P}^3$ and $D \sim 3H$, (so that $D-3H$ is not big and nef) and a singular curve $C\in |D|$ having just one 
singular point $P$ (either a node or a cusp). Then by using Theorem \ref{thm2}  (cfr. examples in \S \ref{severivariety}), $C$ is geometrically linearly normal. We want to show that $C$ is not geometrically $2$-normal. To do this we consider the following chain of implications: 

\begin{eqnarray*}
\mathcal{I}_{N}(D+K_S-H)\,\,\,\,\,  0-\textrm{regular}& 
\Rightarrow\mu_{0,\tilde{C}}\,\,\,\,\textrm{surjective}\Rightarrow\,\,\,\,\rho (g,r,d)\leq 0 & \\
\Downarrow \,\,\,\,\,\,\,\,\,\,\,\,\,\,\,\,\,\,\,\,\,\,\,\,\,\,\,\,\,\,\,\,\,\,& &\\
\,\,\,\,\,\,h^1(\mathcal{I}_{N}(D+K_S-2H))=0& &\\
\Updownarrow\,\,\,\,\,\,\,\,\,\,\,\,\,\,\,\,\,\,\,\,\,\,\,\,\,\,\,\,\,\,\,\,\,\,& &\\
C\,\,\,\,\,\textrm{geom.}\,\, 2-\textrm{normal}\,\,\,\,\,& &\\ 
\end{eqnarray*}
In particular, observe that the 0-regularity of $\mathcal{I}_{N}(D+K_S-H)$ 
implies, but is not equivalent to the vanishing of $h^1(\mathcal{I}_{N}(D+K_S-2H))$, whereas it would be equivalent if $D-3H$ were big and nef.
Now in this example an immediate computation shows that 
$\rho(g,r,deg(C))>0$, so that $\mathcal{I}_{N}(D+K_S-H)$ cannot be 0-regular. On the other hand the 0-regularity could fail because $h^2(\mathcal{I}_{N/S}(D+K_S-3H))\neq 0$, whereas $h^1(\mathcal{I}_{N/S}(D+K_S-2H))$ still vanishes.
But in our case $h^1(\mathcal{I}_{N}(D+K_S-2H))=h^1(\mathcal{I}_{P}(-H))=1$. Thus such a $C$ is not geometrically $2$ normal. 
\end{example}

\section{Examples and application of geometric $k$-normality to Severi varieties}\label{severivariety}

Let us recall that given a Bertini linear system $|D|$ on a smooth surface $S$,
the Severi variety $\mathcal{V}_{\delta,|D|}$ is the locally closed subscheme
of $|D|$, parametrizing irreducible nodal curves in $|D|$, with $\delta$ nodes and no other singularity. 
Saying that $\mathcal{V}_{\delta,|D|}$ is 
{\em regular} at a point $[C]$ means that $\mathcal{V}_{\delta,|D|}$ is smooth and of the expected dimension at $[C]$ 
(i.e. of $dim_{[C]}(\mathcal{V}_{\delta,|D|})=dim(|D|)-\delta$). 
In particular, the regularity  of the Severi variety  
$\mathcal{V}_{\delta,|D|}$ at a point $[C]$, is equivalent to the fact 
that the nodes of $C$ can be indipendently smoothed in $|D|$.   
Now, generalizing the approach of \cite{z1} as it has been done by 
\cite{cs}, it is easy to see that the tangent space $$T_{[C]}|D|
=H^0(\os(D))/H^0(\os),$$ whereas 
$$T_{[C]}(\mathcal{V}_{\delta,|D|})=H^0(\mathcal{I}_{N}(D))/\mathbb{C},$$ where $N$ is the subscheme of the nodes of $C$. 

From these identifications it follows that $C \subset S$ is a regular point $[C]\in \mathcal{V}_{\delta,|D|}$ if 
and only if the subscheme $N$ impose independent conditions on $|D|$ on $S$. Obviously, a sufficient condition for 
this to hold is that 
$h^1(\mathcal{I}_{N}(D))=0$, which is also necessary if $h^1(\os(D))=0$. 

Now, using the results proved in the previous section, 
and choosing smooth non-degenerate surfaces $S\subset \mathbb{P}^r$
with $h^1(\os(kH))=0$ and $K_S \sim kH$, it is possible to relate the regularity of the Severi variety 
$\mathcal{V}_{\delta, |D|}$ at $[C]$
to the fact that $C$ is g.$k$.n. in $\mathbb{P}^r$; indeed this follows
from Theorem \ref{thm1} and the fact that in this case $|D+K_S-kH|=|D|$.
 
In order to study the regularity of 
some Severi varieties on complete intersection surfaces, we can restrict 
Corollary \ref{exampleprop} to the case of nodal curves.

\begin{example}
We examine again an example given in \cite{flaminio} (page 752-753). Let $S 
\subset \mathbb{P}^3$ be a smooth sextic surface. Then there exists an 
irreducible nodal curve $C\in |8H|$ with 48 nodes as it is proved in 
\cite{flaminio}. Such a curve is geometrically linearly normal as it easily follows from e.g. 
Proposition \ref{exampleprop}, but 
{\em is not} geometrically 2-normal. Indeed, since in this case in the 
exact sequence
\[
0 \rightarrow \mathcal{I}_{N}(D+K_S-2H) \rightarrow \os(D+K_S-2H) 
\rightarrow \mathcal{O}_{N} \rightarrow 0
\]
we have $K_S \sim 2H$, $D \sim 8H$, the vanishing of 
$h^1(\mathcal{I}_{N}(8H))$ is a necessary and sufficient condition for the geometric 2-normality; 
this is due to the fact that $h^1(\os(8H))=0$.   
On the other hand, in the construction of such a curve, one finds 
$h^1(\mathcal{I}_{N}(8H))=1$ (see \cite{flaminio}, page 753), so that $C$ is 
not geometrically 2-normal. Moreover observe that in this case the bound \eqref{exampleprop1} gives exactly 
$\delta <48$, so that the bound given is sharp, at least in this example. 
Thus the corresponding point in $\mathcal{V}_{48,|8H|}$ is not a regular point. 
\end{example}

Recall that, in general, the failure of the regularity of a point $[C]$ of a Severi variety 
can be due to several  geometric situations occuring in $\mathcal{V}_{\delta,|D|}$. 
Either $[C]$ is a singular point of a generically smooth component of 
$\mathcal{V}_{\delta,|D|}$, or $[C]$ belongs to a smooth but superabundant component 
(i.e. of dimension higher than the expected), or even $[C]$ is a point on a 
non-reduced component. This last situation would be the most interesting one, 
since there is no example of a non-reduced component 
of some Severi variety.

To conclude the paper, we consider further interesting remarks on Severi varieties of 
plane curves. Since 
every divisor $D\subset\p$ is a multiple of the hyperplane
section $H\subset\p$, to simplify the notation, we will indicate with
$\v$, instead of 
$\mathcal{V}_{\d, nH}$, the Severi variety of irreducible and reduced plane
curves of degree $n$ with $\d$ nodes as the only singularities.   
In this case, we cannot relate the $k$-geometric normality of a nodal curve
with the regularity of the respective Severi variety at the corresponding 
point. Indeed, $K_{\p} \sim -3 H$, so $k= -3$. 

Anyhow, it is classically known that the Severi varieties $\v$  
are always regular, irreducible and not empty  for any
$\d\leq\g$, (see \cite{h}, \cite{severi}, Anhang F and \cite{z1}). 
However, observe that, even when $\d\leq\g$, it is not true in general that 
for any chosen $\delta$ points in the plane, there exists an irreducible plane curve 
of degree $n$ with nodes at these points. 
The classic examples are sextics with nine nodes 
in general position\footnote{The only sextic with nine nodes $P_1,\dots, P_9$ in general position 
in the plane is the double cubic passing through $P_1,\dots, P_9$} 
(see \cite{alb}).

We wonder whether there exists, 
for any $k\leq n-3$, an irreducible curve of degree $n$ with $\delta$ nodes 
in sufficiently general position which impose independent linear 
conditions on the linear system of plane curves of degree $n-k-3$. 
The following theorem gives a partial answer to this question and extends some 
results of \cite{ser}.
\begin{theorem}\label{sv}
Let $\v\subset\mathbb P (H^0(\p, \mathcal{O}_{\p}(n))):=\n$, 
with $N=\frac{n(n+3)}{2}$, the Severi variety of plane curves of degree $n$ with $\d$ nodes. 
For all non negative integers $\,\,n$, $\d$ and $k=1,2,3$ such that $n-3-k\geq 0$ and
\begin{equation}\label{n} 
\d  \leq  \frac{n^2-(3+2k)n+2+k^2+3k}{2}=h^0(\p, \,\mathcal{O}_{\p}(n-k-3)),
\end{equation} 
the {\em general} element of $\,\v$ parametrizes a geometrically $k$-normal plane curve. 
Equivalently, if $n$, $\d$ and $k$ are as before, the general element of $\v$
parametrizes a plane curve with nodes in sufficiently general position
to impose independent linear conditions on the system of plane curves of degree $n-3-k$.
\end{theorem}

Before proving the above result, we observe that, 
from Corollary \ref{exampleprop}, if we have $\d <\frac{n(n-2k)}{4}$, 
then Theorem \ref{sv} is true for any $k$ such that $n-3-k\geq 0$ and for any element in $\v$. 

In the case $k=1$, the Theorem \ref{sv} has already been proved by Sernesi in \cite{ser}. 
Our proof generalizes his result. 
As we shall see, the induction on the 
degree of the curve is, mutata mutandis, essentially the same of Sernesi, 
but our induction on the number of nodes 
is different. Indeed, we use some results about the
local geometry of the closure of $\v$ in $|\mathcal{O}_{\mathbb{P}^2}(n)|$ 
for which we refer to \cite{severi}, \cite{alb}, and \cite{z1}.
Finally, let us remark that Theorem \ref{thm1} implies that 
the bound given in (\ref{n}) cannot be improved. 

\begin{proof}
First of all, let us observe that, if $k$ is an arbitrary
positive integer such that \,\,\,\,$n-3-k\geq 0$ and, if  there exists a point $[C]\in \v$ which 
parametrizes  a geometrically $k$-normal curve with only $\d$
nodes as singularities, then also the general element
of $\v$ parametrizes a g. $k$. n. curve.
To see this, let $\Delta\subset\v$ be a general curve through  $[C]$.
Consider a local analytic parametrization of $\Delta$ at $[C]$, which we will still indicate with $\Delta$.
Taking the restriction to $\Delta$ of the tautological family 
$$
 \{(P, [C]) \mid P\in C\}\subset\p\times\v,
$$
we obtain a family of irreducible plane curves with $\d$ nodes
$$
\psi :\mathcal{C}\to\Delta
$$
parametrized by a smooth curve, whose special fibre $\mathcal{C}_0$ is equal 
to $C$. 
By normalizing $\mathcal{C}$ we obtain a 
family of smooth curves 
$$
\tilde\psi :\tilde{\mathcal{C}}\to\Delta
$$
of geometric genus $g=p_a(C)-\delta$, 
\begin{eqnarray*}
\tilde{\mathcal{C}}\,\,\,\,\,\, &\to &\,\,\,\,\,\,\mathcal{C}\subset\p\times\Delta\\
\searrow & & \swarrow\,\,\,\,\,\, \,\,\, \\
& \Delta & 
\end{eqnarray*}
whose fibres are the normalizations of the curves of the family $\mathcal{C}$.
If $H_z$ is the pullback to $\tilde{\mathcal C}_z :={\tilde\phi}^{-1} (z)$ of a generic line $H\subset\p$, 
using semicontinuity of the cohomology, we have
\begin{eqnarray*}
h^o(\tilde{\mathcal C}_z, \mathcal{O}_{\tilde{\mathcal C}_z}(kH_z)) & \leq & 
h^o(\tilde{\mathcal C}_0,\mathcal{O}_{\tilde{\mathcal C}_0}(kH_0)\\
& = & H^0(\p, \mathcal{O}_{\p}(k )).
\end{eqnarray*}
Then, so $C_z:={\phi}^{-1} (z)$ is geometrically $k$-normal.
In the same way, if the claim is true for fixed $n$, $k\leq n-3$,  
and $\d$ as in \eqref{n}, then the claim
is true also for any ${\d }^\prime\leq \d$.
In fact, as Severi showed in \cite{severi}, if $\d^\prime\leq\d\leq\g$,
then $\v\subset\overline{V_{n,\d^\prime}}$. In particular, if 
$[C]\in\v$ parametrizes a curve with nodes $P_1, ..., P_{\d}$, then there are  
${\d \choose {\d}^{\prime}}$ analytic smooth branches of
${V_{n,\d^\prime}}$ passing through $[C]$ and intersecting transversally at $[C]$.  
Each of these branches parametrizes a family 
of plane curves tending to $[C]$  with $\d^\prime$ nodes, specializing 
to a given subset of $\d^\prime$ nodes of $C$. 
Then, let $\Delta$ be a curve through $[C]$ contained in one of the  
branches of ${V_{n,\d^\prime}}$ and, as before, let us consider the tautological family of projective curves 
$$
\phi :\mathcal{C}\to\Delta
$$
parametrized by $\Delta$. If we normalize $\mathcal{C}$ we obtain a family of curves
$$
\tilde\phi :\tilde{\mathcal{C}}\to\Delta
$$
whose generic fibre is the normalization of the generic fibre
of $\mathcal{C}$, and whose special fibre $\tilde{\mathcal C}_0$ is a partial normalization of the original curve $C$.  
Due to the geometric $k$-normality of $C$, we have  $h^0({\tilde{\mathcal C}}_0,\mathcal{O}_{\tilde{\mathcal C}_0}(kH))=\kn$, with the obvious notation.
Applying semicontinuity to $\tilde{\mathcal{C}}$, we have that
also the generic fibre of $\mathcal{C}$ is geometrically $k$-normal.
Finally, to prove the theorem, it's enough to prove the claim when the equality holds in (\ref{n}). Furthermore, from now on, we shall suppose $k=1,\,2$ or $3$.

Suppose the theorem true for $n$ and let $[\Gamma]\in\v$ correspond to a geometrically $k$-normal curve. Then, let $D$ be a smooth plane curve of degree $k$ which tranversally intersects $\Gamma$ 
and let $P_1, ...,P_{\frac{k^2+3k+2}{2}}$ be $\frac{k^2+3k+2}{2}$ marked
points of $\Gamma\cap D$. If $\Gamma^\prime =\Gamma\cup D\subset\p$, then 
$P_1, ...,P_{\frac{k^2+3k+2}{2}}$ 
are nodes for $\Gamma^\prime$. Let $C\to \Gamma$ be the normalization
of $\Gamma$ and $\cp \to \Gamma^\prime$ the partial normalization of 
$\Gamma^\prime $, obtained by smoothing all singular points of $\Gamma^\prime$, except $P_1, ...,P_{\frac{k^2+3k+2}{2}}$. 
We have the following exact sequence of sheaves on $\cp$,
\begin{equation}\label{se}
0\to\mathcal O_{D}(kH)(-P_1- ...-P_{\frac{k^2+3k+2}{2}})\to
\mathcal O_{\cp}(kH)\to\mathcal O_{C}(kH)\to 0,
\end{equation}
where $kH$ is the pullback by $\cp \to \Gamma^\prime$ of the generic
line in $\p$. Since $1\leq k\leq 3$, then 
$$
deg(\mathcal O_D(kH)(-P_1- ...-P_{\frac{k^2+3k+2}{2}}))<0.
$$ 
Thus,
$$
h^0(D,\mathcal O_D(kH)(-P_1- ...-P_{\frac{k^2+3k+2}{2}}))=0
$$
and so
\begin{eqnarray*}
h^0(\cp, \mathcal O_{\cp}(kH)) & \leq & h^0(C, \mathcal O_{C}(kH)\\
& = & h^0(\p,\mathcal O_{\p}(kH)).
\end{eqnarray*}
But, naturally, $h^0(\cp, \mathcal O_{\cp}(kH)) \geq  h^0(\p, \mathcal O_{\p}(kH))$.\\ Thus,
$h^0(\cp, \mathcal O_{\cp}(kH)) = h^0(\p, \mathcal O_{\p}(kH))$. 
Now, it follows from what Severi proved in \cite{severi} that 
$\Gamma^\prime\in \overline{V}_{n+k, \d+nk-\frac{k^2+3k+2}{2}}$.
In particular, we can obtain $\Gamma^\prime$ as the limit of a family of 
irreducible plane curves 
$$
\psi : \mathcal C \to \Delta
$$
of degree $n+k$ with 
$
\d +nk-\frac{k^2+3k+2}{2}=\frac{(n+k)^2-(3+2k)(n+k)+k^2+3k+2}{2}
$ 
nodes specializing to nodes of $\Gamma^\prime$ different from the marked points $P_1, ...,P_{\frac{k^2+3k+2}{2}}$.
Normalizing $\mathcal C$, we 
obtain a family whose special fibre is exactly $\cp$, and we conclude the induction on the degree via semicontinuity as before.

Now, we have only to show the first step of the induction. Assuming always $k=1, 2$ or $3$, at the first step of induction, 
we have the pairs of values $(n, \d)=(k+a+3, \frac{k^2+3k+2}{2})$, with $a=0, ..., k-1$. For $a=0$, 
by using theorem \ref{thm1}, the claim is true because one point imposes independent conditions on
regular functions. Suppose that the theorem is true for $a<k-1$ and $k=2, 3$, and let $\Gamma$ be a g.$k$.n. plane curve of degree $k+a+3$ with $\frac{a^2+3a+2}{2}$ nodes. 
If $D$ is a line which  intersects transversally $\Gamma$, and if we choose $k+1$ points $P_1, ..., P_{k+1}$ of intersection 
between $\Gamma$ and $D$, we have that the cohomology group $H^0(D,\mathcal O_D(k)(-P_1, ..., -P_{k+1}))$ is zero.
Now, let $C\to\Gamma$ be the normalization of $\Gamma$,  
$C^\prime\to\Gamma^{\prime}$ the partial normalization of $\Gamma^\prime =\Gamma\cup D$ obtained  
smoothing the nodes of $\Gamma\cup D$ different from $P_1, ..., P_{k+1}$, and 
$\mathcal C\to\Delta$ a family of plane curves of degree $k+a+4$
with $\frac{a^2+3a+2}{2}+a+2=\frac{(a+1)^2+3(a+1)+2}{2}$ nodes as only
singularities, tending to $\Gamma\cup D$ in such a way that the nodes of the
general curve of the family tend to the nodes of $\Gamma^\prime$ different 
to the marked $k+1$ points. By using an exact sequence like (\ref{se}),  
$H^0(C^\prime, \mathcal O_{C^\prime}(kH)=H^0(C, \mathcal O_{C}(kH)$.
Then, we conclude as before, applying the semicontinuity theorem to the family
of curves obtained normalizing the family $\mathcal C\to\Delta$. 
\end{proof}

We are persuaded that this theorem admits a generalization to algebraic systems
of plane curves with nodes and a limited number of cusps. But this requires a 
deeper study of the boundary of these varieties, which we shall not consider 
in this paper.

\vskip 5pt

\noindent
{\bf Acknowledgements}: It is our pleasure to thank Dr. Flaminio Flamini for
having suggested the problem and for extremely valuable
discussions and suggestions on these topics. We also wish to thank Proff. I. Dolgachev 
and A. Verra, for the fine lectures they gave during the Pragmatic 2003 school 
activities, and the local organisers (especially 
Prof. Ragusa), for the nice atmosphere.

\end{document}